\documentclass[11pt]{amsart}

\usepackage[margin=1.05in]{geometry}
\usepackage{amsmath,amssymb,amsthm,mathtools, tikz, verbatim}
\usepackage{enumitem}
\usepackage{hyperref}

\newtheorem{theorem}{Theorem}[section]
\newtheorem*{theorem*}{Theorem}

\newtheorem{lemma}[theorem]{Lemma}
\newtheorem{corollary}[theorem]{Corollary}
\newtheorem{remark}[theorem]{Remark}

\newcommand{\R}{\mathbb R}
\newcommand{\N}{\mathbb N}
\newcommand{\Z}{\mathbb Z}
\newcommand{\Q}{\mathbb Q}
\newcommand{\W}{\mathcal W}
\newcommand{\Bad}{\mathbf{Bad}}
\newcommand{\Ccal}{\mathcal C}
\newcommand{\Scal}{\mathcal S}
\newcommand{\Tcal}{\mathcal T}
\newcommand{\Pcal}{\mathcal P}
\newcommand{\Rcal}{\mathcal R}
\newcommand{\dimH}{\dim_{\mathrm H}}
\newcommand{\dimP}{\dim_{\mathrm P}}
\newcommand{\dist}{\operatorname{dist}}

\newcommand{\Leb}{\lambda}
\newcommand{\bp}{\mathbf{p}}
\newcommand{\cL}{\mathcal{L}}
\newcommand{\cV}{\mathcal{V}}

\title{Cartesian product of $\psi$-badly approximable sets}

\author[Mumtaz Hussain]{Mumtaz Hussain}
\address{Mumtaz Hussain,  Department of Mathematical and Physical Sciences,  La Trobe University, Bendigo 3552, Australia. }
\email{m.hussain@latrobe.edu.au}

\author[Johannes Schleischitz]{Johannes Schleischitz}
\address{Johannes Schleischitz, School of Computer, Data, and Mathematical Sciences, Western Sydney University, NSW, Australia}
\email{J.Schleischitz@westernsydney.edu.au}

\author[B. Ward]{ Benjamin Ward}
\address{Benjamin Ward,  {Department of Mathematics, University of York, Deramore Lane, York, YO10 5GH, UK.} }
\email{benjamin.ward@york.ac.uk, ward.ben1994@gmail.com }

\begin{document}
\maketitle

\begin{abstract}
In this note, we study the Cartesian products of sets in $\R^{d}$ that are $\psi$-approximable but not $c\psi$-approximable for arbitrarily small $c>0$, often denoted $\Bad_{d}(\psi)$. Setting $d=\sum_{i=1}^{N}d_i\in\N$ and $\psi_i(q)=q^{-\tau_i}$ with $\tau_i\geq 1/d_i$, $i=1,\ldots,N$, we calculate
\[
\dim_{\mathrm H}\prod_{i=1}^{N}\mathbf{Bad}_{d_i}(\psi_i)
=\min_{1\leq i\leq N}\left\{\frac{d_i+1}{\tau_i+1}+(d-d_i)\right\}=d- \max_{1\leq i\leq N}\left\{ d_i-\frac{d_i+1}{\tau_i+1} \right\},
\]
for $\dim_{\mathrm H}$ denoting the Hausdorff dimension. The lower bound is obtained by combining the delayed-pruning construction of Koivusalo-Levesley-Ward-Zhang (IMRN, 2024)
with the product Cantor-set method of Wang-Wu (TAMS, 2024). 

\end{abstract}


\section{Introduction and statement of results} \label{introd}

In 1962 Erd\H{o}s \cite{Erdos1962} proved that the sumset $\cL +\cL$, where $\cL$ denotes the set of Liouville numbers
\begin{equation*}
    \cL:=\bigcap_{N=1}^{\infty}\left\{x \in \R: \left| x-\frac{p}{q} \right|<q^{-N} \text{  for infinitely many } \, \, \frac{p}{q}\in \Q\right\},
\end{equation*}
is the real line. That is, $\cL+\cL=\R$. An easy consequence of this and Tricot's dimension formula for the Cartesian product of sets \cite[Theorem 3]{Tricot1982} tells us that \begin{equation*}
    \dim_{\mathrm H} \cL \times \cL =1.
\end{equation*}  
See for example \cite[Sect. 2]{JohRM23} for a proof of this and a higher dimensional generalisation. This naturally leads one to studying the metric properties of Cartesian products of other sets of Diophantine interest. In particular, \cite[Conjecture 1]{JohRM23} asks for the Hausdorff dimension of the Cartesian product of sets of numbers with fixed irrationality exponents. Recall, the set of numbers $x \in [0,1]^{k}$ with irrationality exponent $\frac{1}{k}\leq \mu$ is defined as
\begin{equation*}
    \cV_{k}(\mu):=\left\{x \in [0,1]^{k}: \mu=\sup\left\{\alpha>0: \left\| x -\frac{\bp}{q}\right\|_{\infty}<q^{-\alpha-1} \quad \text{ for infinitely many} \, \, \frac{\bp}{q}\in \Q^{k}\right\} \right\}.
\end{equation*}
Here and throughout, we work with sup norm balls, that is,  $\|x\|_\infty=\max_i |x_i|$. Note that by the equivalence of norms in $\R^{k}$ the set $\cV_k(\mu)$ is independent of norm. Related to $\cV_{k}(\mu)$ we have the set of $\tau$-well approximable numbers.
For $\tau\geq \frac{1}{k}$
define
\begin{equation*}
\W_k(\tau)
 :=\left\{x\in[0,1]^k:
 \left\|x-\frac pq\right\|_\infty<q^{-\tau-1}
 \text{ for infinitely many }(p,q)\in\Z^k\times\N\right\}.
\end{equation*}
It is easy to see that $W_{k}(\tau)\supset \cV_{k}(\mu)$ for every $\mu>\tau\geq \frac{1}{k}$.
Metrical properties of these sets are well-known. For instance, the Jarn\'ik-Besicovitch theorem states that, for $\tau\geq \tfrac{1}{d}$,
\begin{equation*}
    \dimH \W_{d}(\tau)=\frac{d+1}{\tau+1}.
\end{equation*}
Motivated by the work of \cite{Erdos1962} it was proven by Wang and Wu \cite[Corollary 1]{WangWuTAMS}, for $d_{i}\in \N$, $\tau_{i}\geq \tfrac{1}{d_{i}}$, with $i=1,...,N$ and $d=\sum_{i}d_{i}$, that
\begin{equation} \label{Wangwu dim}
    \dim_{\mathrm H} \left( \prod_{i=1}^{N}\W_{d_{i}}(\tau_{i}) \right) = \min_{1\leq i\leq N}\left\{\frac{d_i+1}{\tau_i+1}+(d-d_i)\right\}
    =d- \max_{1\leq i\leq N}\left\{ d_i-\frac{d_i+1}{\tau_i+1} \right\}.
\end{equation}
 In particular, this verifies the upper bound of \cite[Conjecture 1]{JohRM23}. Since the work of \cite{WangWuTAMS} there has been a rapid production of papers on the metric properties of the Cartesian products of various sets \cite{Cheng25, ChengZhang27,He26,GHLS}. The motivation of this note is to show that the same result holds true for the smaller set of $\tau$-badly approximable points defined as follows. For $\tau\geq \tfrac{1}{k}$ define
\begin{equation*}
\Bad_k(\tau)
 :=\W_k(\tau)\setminus\bigcap_{c=1}^{\infty}
       \W_k(q\mapsto c^{-1}q^{-\tau}).
\end{equation*}
That is, $x\in \Bad_k(\tau)$ if and only if $x\in\W_k(\tau)$ and
there exists $c(x)>0$ such that
\begin{equation*}
\left\|x-\frac pq\right\|_\infty\geq c(x)q^{-1-\tau}
\qquad\text{for every }(p,q)\in\Z^k\times\N.
\end{equation*}
For background on this problem and its historical development,  we refer the reader to \cite{KLWZ}. Note that for every $\infty>\mu>\tau\geq \frac{1}{k}$.
\begin{equation} \label{eq: bounds between sets}
    \W_{k}(\tau)\supset \cV_{k}(\mu) \supset \Bad_{k}(\mu).
\end{equation}
Our main result is the following:




\begin{theorem}\label{thm:finite-products} \label{thm:main}
Let $N\geq2$, let $d_1,\ldots,d_N\in\N$, and put
\[
 d=d_1+\cdots+d_N.
\]
Let $\tau_i\geq 1/d_i$ for $1\leq i\leq N$.  Then
\[
\dimH\prod_{i=1}^N\Bad_{d_i}(\tau_i)
=\min_{1\leq i\leq N}\left\{
\frac{d_i+1}{\tau_i+1}+(d-d_i)\right\}
=d- \max_{1\leq i\leq N}\left\{ d_i-\frac{d_i+1}{\tau_i+1} \right\}.
\]
\end{theorem}

\begin{remark}\rm \label{rm: cart exact}
 The particular case of $d_i=1$ and $\tau_i>1$ for all $1\leq i \leq N$ has already been proven in the recent preprint \cite{GHLS}. The method in~\cite{GHLS} is different from ours and is specific to $d_i=1$. In particular, it crucially uses that the rational points $\Q$ are a \textit{well-distributed system} in $\R$ (see for instance \cite[Definition 1.3]{BandiGhoshNandi}). This is not true in higher dimensions.
\end{remark}

As an easy corollary we can deduce the following statement also involving packing dimension, denoted by $\dim_P$.

\begin{corollary} \label{packing dim corollary}
    Let $N, d_{1}, \ldots d_{N}$ and $\tau_{1}, \ldots, \tau_{N}$ satisfy the same properties as in the theorem. Then
    \begin{enumerate}
        \item[i)] $\dim_P  \prod_{i=1}^{N} \Bad_{d_i}(\tau_i)  = d$,
        \item[ii)] $\dim_H \prod_{i=1}^{N}\cV_{d_{i}}(\tau_{i})=d-\max_{1\leq i \leq N}\left\{ d_{i}-\frac{d_{i}+1}{\tau_{i}+1} \right\}$.
    \end{enumerate}
\end{corollary}

\begin{remark}\rm
   For $N=1$ an alternative proof of statement $i)$ was given in \cite[Theorem 2.1]{Joh2025MZ} using the variational principle. See also \cite[Theorem 2.3]{Joh2025MZ} for the packing dimension of the smaller sets of exact approximation order under some restrictions. 
\end{remark}

\begin{remark} \rm
Statement $ii)$ resolves \cite[Conjecture 1]{JohRM23} in full generality. In particular, it shows that the upper bound in equation (12) of \cite{JohRM23} is the correct dimension. Note by Remark~\ref{rm: cart exact} the one dimensional version ($d_1=d_2=\ldots=d_N=1$) is resolved in \cite{GHLS} by the same argument as in the proof of Corollary~\ref{packing dim corollary} below. The case $N=1$ of $ii)$ was proven by G{\"u}ting \cite{Guting63}. That is, it was shown that $\dimH \cV_{d}(\tau)=\dimH \W_{d}(\tau)$. Since then sharper notions of irrationality exponents have been introduced, including sets of logarithmic exact \cite{BDV01}, exact \cite{Bugeaud03,Bugeaud08,BugeaudMoreira11, Fregoli24, BandideSaxce25}, and quantitative exact \cite{BakerWard2025} approximation order. 
\end{remark}

\begin{proof}[Proof of Corollary~\ref{packing dim corollary}]
The packing dimension result follows by applying the Hausdorff
dimension result in an $N+1$ product with the 
additional set being $\Bad_{d_{N+1}}(\tau_{N+1})$ with $\tau_{N+1}$ and $d_{N+1}$ chosen such that
\begin{equation} \label{eq: smallest tau}
    \frac{d_{N+1}+1}{\tau_{N+1}+1}-d_{N+1}< \min_{1\leq i\leq N}\left\{
\frac{d_i+1}{\tau_i+1}-d_i\right\}.
\end{equation}
This can be obtained by taking $d_{N+1}=\max_{1\le i\le N} d_i + 1$ and $\tau_{N+1}$ large enough for the left hand side to be smaller than $\min_{1\le i\le N} (-d_i)=-\max_{1\le i\le N} d_i$. Precisely, by an estimate of Tricot \cite{Tricot1982} (see Theorem~\ref{thm:MT} below) we have that
\begin{align*}
    \dim_P \prod_{i=1}^{N}\Bad_{d_{i}}(\tau_{i}) &\geq \dimH \prod_{i=1}^{N+1}\Bad_{d_{i}}(\tau_{i}) - \dimH \Bad_{d_{N+1}}(\tau_{N+1})\\
    &= \sum_{i=1}^{N+1}d_{i}+\min_{1\leq i\leq N+1}\left\{
\frac{d_i+1}{\tau_i+1}-d_i\right\} -\frac{d_{N+1}+1}{\tau_{N+1}+1}\\
&= \sum_{i=1}^{N}d_{i}=d,
\end{align*}
where the penultimate equation uses Theorem~\ref{thm:main} and the last line uses \eqref{eq: smallest tau}.  

For the Hausdorff dimension statement take any $\varepsilon>0$ and use \eqref{eq: bounds between sets} to see that
\begin{equation*}
    \dim_H \prod_{i=1}^{N} \W_{d_i}(\tau_{i}-\varepsilon) \geq \dim_H \prod_{i=1}^{N}\cV_{d_{i}}(\tau_{i}) \geq \dim_{H}\prod_{i=1}^{N}\Bad_{d_{i}}(\tau_{i}).
\end{equation*}
The lower bound immediately follows from Theorem~\ref{thm:finite-products}. The upper bound follows from \cite[Corollary 1]{WangWuTAMS} and the fact that $\varepsilon>0$ can be chosen arbitrarily small.
\end{proof}

\medskip

The {\em upper bound} of our main result follows directly from the result of \cite{WangWuTAMS}, since it can be seen that $\psi$-badly approximable points are a subset of $\psi$-well approximable points. 
For the {\em lower bound} of Theorem~\ref{thm:main}, which carries the main weight of the paper, we construct a
Cantor subset inside
 $\prod_{i}\Bad_{d_i}(\tau_i)$ inspired by a combination of ideas from \cite{KLWZ,WangWuTAMS}. Briefly, the idea is as follows:
\begin{enumerate}[label=\textup{(\arabic*)},leftmargin=2.3em]
\item Begin with some (sup-norm) ball $B=\prod_{i=1}^{N}B_{i} \subset \R^{d}$. Each of the coordinate spaces take turns in being active and passive, in a similar spirit to \cite{WangWuTAMS}, as described below:

\item \emph{Active coordinate:} Suppose $1\leq i \leq N$ is the active coordinate. In the active coordinate we start with the ball $B_{i}\subset \R^{d_{i}}$ and use the same argument as in \cite{KLWZ} to construct a collection of balls that avoid $c\psi(q)$-neighbourhoods of all rational points $\tfrac{\bp}{q}\in \Q^{d_{i}}$ up to some height on the denominator, while simultaneous being contained in a $\psi(q)$-neighbourhood of a finite number of rational points $\tfrac{\bp}{q} \in \Q^{d_{i}}$. Say $B_{i,a}$ is a ball in this collection $\mathcal{B}_{i}$. We now have a collection of hyperrectangles of the form $\prod_{j=1}^{i-1}B_{j} \times B_{i,a}\times \prod_{j=i+1}^{N}B_{j}$.

\item \emph{Passive coordinate:} In each of the $1\leq j\neq i\leq N$ passive coordinates partition $B_{j}\subset \R^{d_{j}}$ up into a disjoint collection of balls, say $\mathcal{B}_{j,a}$ with radius comparable to $B_{i,a}$. The collection should cover a positive proportion (in terms of Lebesgue measure) of the $B_{j}$ coordinate axis. We have to be careful to ensure the subcollection of balls we choose from this set avoids $c\psi(q)$-neighbourhoods of rational points $\tfrac{\bp}{q}\in \Q^{d_{j}}$. Using results from \cite{KLWZ} we show it is possible to keep a collection of balls with total Lebesgue measure comparable to the Lebesgue measure of $B_{j}$. The collection of approximate balls
\begin{equation*}
    \bigcup_{B_{i,a}\in \mathcal{B}_{i}}
    \ \bigcup_{\substack{B_{j,a,b_j}\in\mathcal{B}_{j,a}\\j\neq i}}
    \left(\prod_{j=1}^{i-1}B_{j,a,b_j}\right)\times B_{i,a}
    \times\left(\prod_{j=i+1}^{N}B_{j,a,b_j}\right)\subset B,
\end{equation*}
forms the first level of our Cantor set. 

\item \emph{Switching:} At the end of this level construction
we switch the active coordinate to the $i+1$ set, modulo $N$.  Repeating this process, this is sufficient to build a Cantor subset $K$ of $\prod_{i=1}^{N}\Bad_{d_i}(\tau_i)$.

\end{enumerate}
Constructing a probability measure on $K$ and applying the mass distribution principle obtains our desired dimension lower bound.
\\

\noindent \textbf{Notation:} Throughout $\Leb_{k}$ denotes $k$-dimensional Lebesgue measure, $|\cdot|$ denotes the diameter of a set, $\dimH$ denotes Hausdorff dimension and $\dim_{\mathrm P}$ denotes packing dimension.\\

\noindent \textbf{Acknowledgments:} Part of this work was done at Ergodic Theory, Diophantine Approximation and Related Topics workshop at Creswick in June 2026. We thank the organiser of the workshop and MATRIX Research Institute for their hospitality and support. The third named author is supported by a Leverhulme Trust Early Career Fellowship (ECF-2024-401).\\

\noindent \textbf{AI Declaration:} The principal mathematical ideas and proofs presented in the paper originate with the authors. ChatGPT5.6 Sol has been used to improve language and assisting the authors in checking arguments that they have developed and written themselves. In all cases, the authors have independently verified the final text and take full responsibility for its correctness.

\section{Preliminaries} \label{upb}

Firstly we dispose of the immediate upper bound of Theorem~\ref{thm:main}.

\begin{proof}[Upper bound in Theorem \ref{thm:main}]
We have that, for any $k\in \N$ and $\tau\geq 1/k$, 
\begin{equation*}
    \W_{k}(\tau)\supseteq \Bad_{k}(\tau), 
\end{equation*}
and so
\begin{equation*}
    \dimH \prod_{i=1}^{N} \Bad_{d_{i}}(\tau_{i}) \leq \dimH \prod_{i=1}^{N}\W_{d_{i}}(\tau_{i})\overset{\eqref{Wangwu dim}}{=}\min_{1\leq i\leq N}\left\{
\frac{d_i+1}{\tau_i+1}+(d-d_i)\right\}.
\end{equation*}
\end{proof}
The rest of the paper is devoted to the proof of the corresponding lower bound statement. For the later use put
\[
 s_i=\frac{d_i+1}{\tau_i+1},\qquad
 \Sigma_i=s_i+d-d_i,
 \qquad
 \Sigma=\min_{1\leq i \leq N}\Sigma_i.
\]
Thus Theorem~\ref{thm:main} says 
\begin{equation*}
    \dimH \prod_{i=1}^{N}\Bad_{d_{i}}(\tau_{i})=\Sigma.
\end{equation*}
The following results due to Tricot and Marstrand \cite{Tricot1982, Marstrand1954} helps us to reduce the problem.

\begin{theorem}[{\cite[Theorem 3]{Tricot1982}}]\label{thm:MT}
Let $A\subset\R^{n}$ and $B\subset\R^{m}$ be Borel sets.  Then
\[
\dimH A+\dimH B\leq \dimH(A\times B)
\leq \dimH A+\dimP B.
\]
\end{theorem}
We may assume that 
\begin{equation}\label{eq:strict}
\tau_i>\frac1{d_i}\qquad(i=1,\ldots, N).
\end{equation}
To see why,  suppose for some $1\leq j \leq N$ we have $\tau_{j}=1/d_{j}$. Then $\Bad_{d_{j}}(\tau_{j})=\Bad_{d_{j}}$ the set of classical badly approximable points. It is well-known that $\dimH \Bad_{d_{j}}=\dim_P \Bad_{d_{j}}=d_{j}$, \cite{Schmidt2}. Therefore using Theorem~\ref{thm:MT} we have that
\begin{equation*}
    \dimH \underset{i\neq j}{\prod_{i=1}^{N}}\Bad_{d_{i}}(\tau_{i})+d_{j}= \dimH \prod_{i=1}^{N}\Bad_{d_{i}}(\tau_{i}).
\end{equation*}
Clearly, if we show that Theorem~\ref{thm:main} holds for the set $\underset{i\neq j}{\prod_{i=1}^{N}}\Bad_{d_{i}}(\tau_{i})$ then
\begin{align*}
    \dimH \prod_{i=1}^{N}\Bad_{d_{i}}(\tau_{i})&=\left(d-d_{j}+\min_{1\leq i\neq j\leq N}\left\{ \frac{d_{i}+1}{\tau_{i}+1}-d_{i} \right\}\right) +d_{j},\\
    &=\min_{1\leq i\leq N}\left\{ \frac{d_{i}+1}{\tau_{i}+1}+(d-d_{i}) \right\}
\end{align*}
Repeating this until we have exhausted all cases where $\tau_i=1/d_i$ we see that we are left to prove Theorem~\ref{thm:main} under assumption \eqref{eq:strict}.
\subsection{The one-coordinate construction from \cite{KLWZ}} \label{sec:t-bad construction}
This section is a repetition of the construction in \cite{KLWZ}. For completeness we include it here since Lemmas~\ref{lem:avoid}-\ref{lem:filling} are essential to our proof. Here, $k$ denotes the dimension of a single coordinate space and
$\tau>1/k$.  
Choose an integer $t$ with $t^k>k!$, and then choose a large integer $N$.
At level $n=0$ we regard $[0,1]^k$ as the single cube. Partition $[0,1]^{k}$ up into $t^{k}N^{(k+1)}$ cubes of equal side lengths and let $\Ccal_{1}$ denote this collection of cubes. At the next level take each cube $I \in \Ccal_{1}$ and partition it into $N^{(k+1)}$ cubes of equal side lenghts, let $\Ccal_{2}(I)$ denote this set of cubes and $\Ccal_{2}=\bigcup_{I\in \Ccal_{1}}\Ccal_{2}(I)$. Repeating this process within each cube of the previous level, at level $n\geq1$ the cubes $I \in \Ccal_{n}$ have side lengths
\begin{equation}\label{eq:rho-level}
\rho_n=t^{-1}N^{-n(k+1)/k}.
\end{equation}
\begin{lemma}[Simplex lemma, \cite{KTV}]\label{lem:simplex}
Let $Q\in\N$ and let $E\subset\R^k$ be convex with
\[
\Leb_k(E)<(k!)^{-1}Q^{-(k+1)}.
\]
If $E$ contains $k+1$ rational points with denominators at most $Q$, then
these rational points lie on some hyperplane.
\end{lemma}
Since
\[
\Leb_k(I)=t^{-k}N^{-n(k+1)}<(k!)^{-1}(N^n)^{-(k+1)}
\]
for every $I\in\Ccal_n$, Lemma \ref{lem:simplex} implies that all rational
points $p/q\in I$ with $1\leq q<N^n$ lie on a hyperplane.
Applying this to each cube $I\in \Ccal_{n}$ let $\cL_n$ denote the
collection of all hyperplanes obtained at level $n$.

Set $u=4$ (as appearing in \cite[eq. (1)]{KLWZ}) and define
\begin{equation}\label{eq:ell-delta}
\ell(n)=\left\lfloor (n+3)\frac{(1+\tau)k}{k+1}\right\rfloor,
\qquad
\delta(n)=t^{-1}N^{-(n+3)(1+\tau)}.
\end{equation}
The floor in the definition of $\ell(n)$ gives
\begin{equation}\label{eq:delta-rho}
N^{-(k+1)/k}\rho_{\ell(n)}\leq \delta(n)\leq \rho_{\ell(n)}.
\end{equation}
For a hyperplane $L\in\cL_n$ and $\delta>0$ write
\[
L^{\delta}=\{x\in[0,1]^k:\dist_\infty(x,L)\leq\delta\}.
\]
 At the first pruning level
$\ell(1)$, remove every level $\ell(1)$ cube that meets $L^{\delta(1)}$
for some $L\in\cL_1$.  The remaining cubes form $\Scal_{\ell(1)}$.
Suppose the construction is complete up to $\ell(n-1)$.  At stage $n$ we
consider only those hyperplanes $L\in\cL_n$ for which there is a
rational point $p/q\in L$ satisfying
\begin{equation}\label{eq:Lstar}
N^{n-1}\leq q<N^n,
\qquad
B_\infty(p/q,\delta(n))\cap
\bigcup_{I\in\Scal_{\ell(n-1)}}I\neq\varnothing.
\end{equation}
Call this subcollection $\cL_n^*$.  Then $\Scal_{\ell(n)}$ is the family
of level $\ell(n)$ cubes contained in some cube of $\Scal_{\ell(n-1)}$ and
meeting none of the thickenings $L^{\delta(n)}$ with $L\in\cL_n^*$.
Set
\begin{equation}\label{eq:Ctau}
\Ccal_k^\tau(N)=\bigcap_{n=1}^{\infty}
     \bigcup_{I\in\Scal_{\ell(n)}}I.
\end{equation}

The following statements provide us with essential properties of $\Ccal_{k}^{\tau}(N)$ that will be used in our product construction. 
\begin{lemma}[cf. {\cite[Lemma 2]{KLWZ}}]\label{lem:avoid}
Let
\[
c_N=t^{-1}N^{-4(1+\tau)}.
\]
Then
\[
\Ccal_k^\tau(N)\cap
\bigcup_{p/q\in\Q^k\cap[0,1]^k}
B_\infty\bigl(p/q,c_Nq^{-1-\tau}\bigr)=\varnothing.
\]
\end{lemma}

\begin{lemma}[{cf. \cite[Lemma 4]{KLWZ}}]\label{lem:density}
For $N$ sufficiently large there is $0<\varepsilon_N<1$, with
$\varepsilon_N\to0$ as $N\to\infty$, such that every surviving cube
{$I\in\bigcup_r\Scal_{\ell(r)}$} satisfies
\[
\Leb_k(I\cap\Ccal_k^\tau(N))\geq (1-\varepsilon_N)\Leb_k(I).
\]
\end{lemma}

Next, define the rational subsystem
\begin{equation}\label{eq:QNtau}
Q(N,\tau)=\left\{\frac{\bp}{q}\in\Q^k\cap[0,1]^k:
\begin{array}{l}
\text{for some }n\text{ and }L\in\cL_n^*,\ \bp/q\in L,\\
N^{n-1}\leq q<N^n,\\
B_\infty(\bp/q,\delta(n))\cap
\displaystyle\bigcup_{I\in\Scal_{\ell(n-1)}}I\neq\varnothing
\end{array}\right\}.
\end{equation}
For $\bp/q\in Q(N,\tau)$ define the large and small balls
\begin{equation}\label{eq:large-small}
B^+(\bp/q)=B_\infty(\bp/q,3q^{-1-1/k}),
\qquad
B^-(\bp/q)=B_\infty(\bp/q,q^{-1-\tau}).
\end{equation}
Note that the large balls are at ($3$ times) the Dirichlet approximation scale, whereas the small balls are at the $\tau$ approximation scale.

\begin{lemma}[{cf. \cite[Lemma 5]{KLWZ}}]\label{lem:big-cover}
Let
\[
\W_k(Q(N,\tau),3q^{-1/k})
=\left\{x\in[0,1]^k:x\in B^+(\bp/q)
\text{ for infinitely many }\bp/q\in Q(N,\tau)\right\}.
\]
Then, for $\tau>1/k$, we have
\[
\Ccal_k^\tau(N)\subset \W_k(Q(N,\tau),3q^{-1/k}).
\]
\end{lemma}
Consequently, as pointed out in \cite{KLWZ}, by combining Lemma~\ref{lem:big-cover} with Lemma~\ref{lem:density} we see that for every surviving cube {$I\in\bigcup_r\Scal_{\ell(r)}$}
\begin{equation}\label{eq:limsup-density}
\Leb_k\left(I\cap\limsup_{\bp/q\in Q(N,\tau)}B^+(p/q)\right)
\geq (1-\varepsilon_N)\Leb_k(I).
\end{equation}

\begin{lemma}[{cf. \cite[Lemma 6]{KLWZ}}]\label{lem:small-cube}
For every $p/q\in Q(N,\tau)$ there is a surviving cube
$J=J(\bp/q)\in\bigcup_r\Scal_{\ell(r)}$ such that
\[
J\subset B^-(\bp/q)
\]
and
\begin{equation}\label{eq:small-comparable}
\Leb_k(B^-(\bp/q))\leq C_1\Leb_k(J),
\end{equation}
where $C_1$ is independent of $\bp/q$.  
\end{lemma}
Clearly this implies there exists some constant $C_{2}>0$ independent of $\bp/q$ such that
\begin{equation}\label{eq:side-comparable}
C_2^{-1}q^{-1-\tau}\leq |J|
\leq C_2q^{-1-\tau}.
\end{equation}

The following lemma, appearing in \cite{KLWZ}, provides an analogue of the $K_{G,B}$ Lemma from \cite{BeresnevichVelani} for $\limsup$ sets with positive measure.

\begin{lemma}[\cite{KLWZ}, Lemma 7]\label{lem:filling}
Let $\{B_\nu\}_{\nu\in\N}$ be a sequence of sup-norm balls in $\R^k$ with
$r(B_\nu)\to0$.  Suppose that there is $C_0>0$ such that
\[
\Leb_k\left(I\cap\limsup_{\nu\to\infty}B_\nu\right)
\geq C_0\Leb_k(I)
\]
for every surviving cube $I\in\bigcup_r\Scal_{\ell(r)}$.  Then for every such
$I$ and every $G\geq1$ there is a finite disjoint subcollection
$\Tcal_{G,I}\subset\{B_\nu:\nu\geq G\}$ such that every ball of
$\Tcal_{G,I}$ is contained in $I$ and
\[
\Leb_k\left(\bigcup_{B\in\Tcal_{G,I}}B\right)
\geq \kappa\Leb_k(I),
\]
where $\kappa>0$ depends only on $k$ and $C_0$.
\end{lemma}

Note that \eqref{eq:limsup-density} ensures Lemma~\ref{lem:filling} is applicable to the pair $\W_k(Q(N,\tau),3q^{-1/k})$ and $\bigcup_r\Scal_{\ell(r)}$.

As outlined in the sketch of the proof, in the \textit{Passive coordinate} step we need to show that inside a surviving cube at any later level of the surviving cubes a positive proportion of the measure is retained. This essentially follows from Lemma~\ref{lem:density}, but for easy recall we provide the following lemma.
\begin{lemma}[Passive splitting]\label{lem:passive}
Let $I\in\bigcup_m\Scal_{\ell(m)}$ be a surviving cube.  There exist constants
$a_0\in(0,1)$ and $\eta_0>0$ such that for every sufficiently small
$0<r<|I|$ there is a finite family
$\Pcal(I,r)$ of $3r$-disjoint surviving cubes $J\subset I$ with
\[
a_0r\leq |J|\leq r
\]
and
\[
\sum_{J\in\Pcal(I,r)}\Leb_k(J)\geq \eta_0\Leb_k(I).
\]
\end{lemma}
\begin{proof}
Choose $m$ with $\rho_{\ell(m)}\leq r<\rho_{\ell(m-1)}$, where $\rho_{\ell(m)}$ is as in
\eqref{eq:rho-level}.  Note
\begin{equation*}
    N^{-(1+\tau)+\frac{k+1}{k}} \geq \frac{\rho_{\ell(m)}}{\rho_{\ell(m-1)}} \geq N^{-(1+\tau)-\frac{k+1}{k}}.
\end{equation*}
 Therefore, a level $\ell(m)$ cube has side length between
$a_0r$ and $r$ with $a_0=N^{-(1+\tau)-(k+1)/k}$.  Let $\Pcal^{*}(I,r)$ be the family of all level $\ell(m)$ surviving cubes contained in $I$. By the Vitali covering lemma there exists a disjoint subcollection $\Pcal(I,r)\subset \{J: J\in \Pcal^{*}(I,r)\}$ such that
\begin{equation*}
    \bigcup_{J\in \Pcal(I,r)} 9J \supseteq \bigcup_{J\in \Pcal^{*}(I,r)} J
\end{equation*}
and $\Pcal(I,r)$ is $3r$-disjoint. By definition
\begin{equation*}
    \bigcup_{J\in \Pcal^{*}(I,r)} J\supseteq I\cap\Ccal_k^\tau(N).
\end{equation*} 
Therefore, using Lemma~\ref{lem:density}, we have
\[
\sum_{J\in\Pcal(I,r)}\Leb_k(J) = \tfrac{1}{9^{k}} \sum_{J\in\Pcal(I,r)}\Leb_k(9J) \geq \tfrac{1}{9^{k}}\Leb_{k}\left(\bigcup_{J\in \Pcal^{*}(I,r)}J\right)
\geq \tfrac{1}{9^{k}}\Leb_k(I\cap\Ccal_k^\tau(N))
\geq \tfrac{1}{9^{k}}(1-\varepsilon_N)\Leb_k(I).
\]
Thus $\eta_0=9^{-k}(1-\varepsilon_N)$ is admissible.
\end{proof}

\section{Cantor construction for the product set}

We now prove the lower bound in Theorem \ref{thm:main} under assumption \eqref{eq:strict}. This is the part of the proof that follows
the product Cantor construction of Wang and Wu \cite{WangWuTAMS}.

For $i=1,\ldots,N$ run the construction of \S~\ref{sec:t-bad construction} with
$k=d_i$, $\tau=\tau_i$, and a sufficiently large integer $N_i$. Write $\Ccal_i=\Ccal_{d_i}^{\tau_i}(N_i)$,
$Q_i=Q(N_i,\tau_i)$, and $\Scal_{\ell_i(r)}^{(i)}$
for the resulting Cantor subset, rational subsystem, and surviving cubes respectively. For $\bp/q\in Q_i$ put
\[
B_i^+(\bp/q)=B_\infty(\bp/q,3q^{-1-1/d_i}),
\qquad
B_i^-(\bp/q)=B_\infty(\bp/q,q^{-1-\tau_i}).
\]
Let $\rho_i(q)=3q^{-1-1/d_i}$ and $r_i(q)=q^{-1-\tau_i}$. Fix any $0<\varepsilon<\min_{1\leq i\leq N}s_i$ and some $C>0$. We shall construct a Cantor set $K\subset\prod_{i=1}^{N}\Bad_{d_i}(\tau_i)$ and a probability
measure $\mu$ on $K$ such that
\begin{equation}\label{eq:frostman}
\mu(B_\infty(x,r))\leq Cr^{\Sigma-\varepsilon}
\qquad(x\in\R^d,\ 0<r<r_0).
\end{equation}
The mass distribution principle, see for example \cite[Proposition 4.2]{Falconer_book2013}, will then give $\dimH K\geq\Sigma-\varepsilon$, and since
$\varepsilon>0$ is arbitrary the desired lower bound follows.

\subsection{Active and passive coordinates}
At each level,  one coordinate block is active and all other coordinate blocks are passive. Write
\[
i(n)=1+((n-1)\bmod N).
\]
Thus the active coordinate cycles through $1,2,\ldots,N$. Each level, whether active or passive, will end with a collection of rectangles, whose projections onto the $i$th coordinate is a collection of surviving cubes from the construction of $\Ccal_i$. Briefly, when coordinate $i$ is active we use Lemma~\ref{lem:filling} to obtain a collection of large
balls $B_i^+$ and then shrink to the small balls $B_i^-$. When coordinate
$j\neq i$ is passive we use Lemma~\ref{lem:passive} to cut the current surviving cube into smaller surviving cubes at the same scale as the active small ball.

\subsection{Nested rectangles}

We construct finite families $\Rcal_n$ of rectangles
\[
R=\prod_{j=1}^{N}I_j(R)\subset\prod_{j=1}^{N}[0,1]^{d_j},
\]
where each $I_i(R)$ is a surviving cube in the $i$th coordinate
construction. The rectangles are approximate cubes: there is a constant
$C\geq1$ such that
\begin{equation}\label{eq:shape}
C^{-1}|I_i(R)|\leq |I_j(R)|\leq C|I_i(R)|
\qquad(1\leq i,j\leq N)
\end{equation}
for every rectangle in the construction. We start with
\[
\Rcal_0=\left\{\prod_{j=1}^{N}[0,1]^{d_j}\right\}.
\]
Assume $\Rcal_{n-1}$ has been constructed and fix $R=\prod_{j=1}^{N}I_j\in\Rcal_{n-1}$, where $I_j$ is a surviving cube in the $j$th construction. Let $i=i(n)$ be the active coordinate and let every $j\neq i$ be passive. Apply Lemma~\ref{lem:filling} to $I_i$ with some large $G(R,n)>0$ to be determined, to see there is a finite disjoint family
\[
\Tcal(R)=\{B_i^+(\bp_a/q_a)\subset I_i\}_a
\]
such that
\begin{equation}\label{eq:large-fill}
\sum_a\Leb_{d_i}(B_i^+(\bp_a/q_a))
\geq\kappa_i\Leb_{d_i}(I_i).
\end{equation}
Shrink each of the balls in $\Tcal(R)$ from $B_i^+(\bp_a/q_a)$ to $B_i^-(\bp_a/q_a)$. Since $\tau_i>1/d_i$, once $G(R,n)$ is large enough we have
\[
3B_i^-(\bp_a/q_a)\subset B_i^+(\bp_a/q_a)\subset I_i.
\]
Apply Lemma~\ref{lem:small-cube} to each $B_i^-(\bp_a/q_a)$ to obtain some surviving cube
\[
J_{i,a}\subset B_i^-(\bp_a/q_a)
\]
with
\begin{equation}\label{eq:active-size}
|J_{i,a}|\asymp q_a^{-1-\tau_i},
\qquad
\Leb_{d_i}(J_{i,a})\asymp q_a^{-(1+\tau_i)d_i}.
\end{equation}

In each passive coordinate $j\neq i$ use Lemma~\ref{lem:passive} to choose a family
$\Pcal_{j,a}=\Pcal(I_j,q_a^{-1-\tau_i})$ of pairwise $3r$-disjoint surviving cubes
$J_{j,a,b_j}\subset I_j$, such that
\begin{equation}\label{eq:passive-size}
|J_{j,a,b_j}|\asymp q_a^{-1-\tau_i},
\qquad
\sum_{b_j}\Leb_{d_j}(J_{j,a,b_j})\geq\eta_j\Leb_{d_j}(I_j),
\end{equation}
where $\eta_j>0$ is the constant from Lemma~\ref{lem:passive} which depends on $\tau_j,N_j$, but is independent of $I_j$ and $a$.
For each $a$ and each independent choice $b=(b_j)_{j\neq i}$ define a child rectangle $R_{a,b}$ of $R$ by taking
$J_{i,a}$ in the active coordinate and $J_{j,a,b_j}$ in each passive coordinate. That is,
\[
R_{a,b}=\left(\prod_{j=1}^{i-1}J_{j,a,b_j}\right)
\times J_{i,a}\times
\left(\prod_{j=i+1}^{N}J_{j,a,b_j}\right).
\]

The children of a fixed parent are pairwise $3r$-disjoint, since different active
indices $a$ lie inside disjoint large active balls, and for a fixed $a$ in at least one of the passive coordinate we are $3r$-disjoint. Note that
\eqref{eq:active-size} and \eqref{eq:passive-size} easily imply \eqref{eq:shape}.
By taking $G(R,n)$ increasing fast enough we ensure that
\begin{equation}\label{eq:fast-shrink}
\max_a\{q_a^{-1-1/d_i}\}\leq 2^{-n}\min_{1\leq j\leq N}|I_j|.
\end{equation}
Hence the diameters of the construction rectangles tend to zero.

Repeating this construction for each $R\in\Rcal_{n-1}$ we obtain our set of level $n$ rectangles
\[
\Rcal_n=\bigcup_{R\in\Rcal_{n-1}}\Rcal_n(R),
\qquad \Rcal_n(R)=\{R_{a,b}:a,b\text{ as above}\}.
\]

Define the Cantor subset by
\begin{equation}\label{eq:product-K}
K=\bigcap_{n=0}^{\infty}\bigcup_{R\in\Rcal_n}R.
\end{equation}

\begin{lemma}\label{lem:membership}
For any choice of sequence $G(R,n)$ tending to infinity,  the Cantor subset $K$ satisfies
\[
K\subset\prod_{i=1}^{N}\Bad_{d_i}(\tau_i).
\]
\end{lemma}

\begin{proof}
Let $x=(x_1,\ldots,x_N)\in K$. For each $i=1,\ldots,N$, the coordinate $x_i$ lies in a
nested sequence of surviving cubes for the $i$th one-coordinate
construction. Hence $x_i\in\Ccal_i$.

By Lemma~\ref{lem:avoid}, there exists $c_i>0$ such that
\[
\left\|x_i-\frac pq\right\|_\infty\geq c_iq^{-1-\tau_i}
\qquad\text{for all }p/q\in\Q^{d_i}\cap[0,1]^{d_i}.
\]
Thus $x_i$ is not $cq^{-1-\tau_i}$-well approximable for sufficiently small $c$.

It remains only to show that $x_i\in\W_{d_i}(\tau_i)$. Coordinate $i$ is
active at infinitely many levels. At each such level the construction
places $x_i$ inside a small ball
$B_i^-(p/q)=B_\infty(p/q,q^{-1-\tau_i})$ for some $p/q\in Q_i$. The
thresholds $G(R,n)$ tend to infinity, so the denominators occurring along
these active levels tend to infinity. Hence
\[
\left\|x_i-\frac pq\right\|_\infty<q^{-1-\tau_i}
\]
for infinitely many rational points $p/q$. Therefore
$x_i\in\W_{d_i}(\tau_i)$, and the lower bound above shows
$x_i\in\Bad_{d_i}(\tau_i)$.
\end{proof}

\subsection{Mass distribution}

Starting with $\mu([0,1]^d)=1$, recursively define the measure on the rectangles. That is, suppose $\mu(R)$ has been defined for each $R\in\Rcal_{n-1}$ and that $i$ is the active coordinate at level $n$. Then for each child $R_{a,b}$ of $R$ define
\begin{equation}\label{eq:recursive-mass}
\mu(R_{a,b})
=\mu(R)\frac{\Leb_{d_i}(B_i^+(\bp_a/q_a))}{\sum_a\Leb_{d_i}(B_i^+(\bp_a/q_a))}
\prod_{j\neq i}\frac{\Leb_{d_j}(J_{j,a,b_j})}{\sum_{b_j}\Leb_{d_j}(J_{j,a,b_j})}.
\end{equation}
The mass of a parent rectangle is preserved over the total mass of the children rectangles, so this gives a
Borel probability measure supported on $K$.
By \eqref{eq:large-fill} and \eqref{eq:passive-size},
\begin{align}\label{eq:denom-bounds}
&\sum_a\Leb_{d_i}(B_i^+(\bp_a/q_a))\geq\kappa_i\Leb_{d_i}(I_i),\\
&\sum_{b_j}\Leb_{d_j}(J_{j,a,b_j})\geq\eta_j\Leb_{d_j}(I_j)
\qquad(j\neq i).\nonumber
\end{align}
Note $\Leb_d(R)=\prod_{j=1}^{N}\Leb_{d_j}(I_j)$.
From \eqref{eq:recursive-mass} and \eqref{eq:denom-bounds} we have
\begin{align}
\mu(R_{a,b})
&\leq  (\kappa_{i}\prod_{j\neq i}\eta_{j})^{-1}\mu(R)
\frac{q_a^{-(1+1/d_i)d_i}}{\Leb_{d_i}(I_i)}
\prod_{j\neq i}\frac{q_a^{-(1+\tau_i)d_j}}{\Leb_{d_j}(I_j)}\nonumber\\
&\leq (\kappa_{i}\prod_{j\neq i}\eta_{j})^{-1}\frac{\mu(R)}{\Leb_d(R)}
\left(q_a^{-1-\tau_i}\right)^{s_i+\sum_{j\neq i}d_j}\nonumber\\
&=(\kappa_{i}\prod_{j\neq i}\eta_{j})^{-1}\frac{\mu(R)}{\Leb_d(R)}
\left(q_a^{-1-\tau_i}\right)^{\Sigma_i}\nonumber\\
&\leq (\kappa_{i}\prod_{j\neq i}\eta_{j})^{-1} C_{2}^{\Sigma_{i}}\frac{\mu(R)}{\Leb_d(R)}|R_{a,b}|^{\Sigma_{i}}.
\label{eq:one-step}
\end{align}

Observe that if $R=[0,1]^d$, i.e. we are in the first level of the construction of $K$, then we immediately obtain from \eqref{eq:one-step} that for every child $R_{a,b}$
\begin{equation*}
\mu(R_{a,b})\leq C'|R_{a,b}|^{\Sigma_i},
\end{equation*}
with $C'=\max_{1\leq i \leq N}\left\{ (\kappa_{i}\prod_{j\neq i}\eta_{j})^{-1} C_{2}^{\Sigma_{i}}\right\}\geq 1$. Note that the value $G(R,n)$ is chosen after the parent $R$ and its mass
$\mu(R)$ have already been defined. For any $\varepsilon>0$ we can choose $G(R,n)$ sufficiently large, depending on $R$ and $\varepsilon$, so that
\begin{equation}\label{eq:threshold}
C'\frac{\mu(R)}{\Leb_d(R)}\leq
\left(q_a^{-1-\tau_i}\right)^{-\varepsilon}
\qquad\text{for every selected }a.
\end{equation}
This is possible because $q_a^{-1-\tau_i}\to0$. Combining \eqref{eq:one-step} with \eqref{eq:threshold} we obtain that for every child $R_{a,b}$ of parent $R$ we have
\begin{equation*}
\mu(R_{a,b})\leq
\left(q_a^{-1-\tau_i}\right)^{\Sigma_i-\varepsilon}
\leq C|R_{a,b}|^{\Sigma-\varepsilon}.
\end{equation*}
Since $\Sigma_i\geq\Sigma$ for every active coordinate we have the following lemma:

\begin{lemma}\label{lem:rectangle-mass}
For every $\varepsilon>0$ sufficiently small there exists a sequence $G(R,n)$ such that every rectangle in the construction of $K$ satisfies
\[
\mu(R)\leq C'|R|^{\Sigma-\varepsilon}.
\]
\end{lemma}

\section{The H\"older exponent of $\mu$}

We need the following general lemma:

\begin{lemma}[{\cite[Lemma 7]{BeresnevichVelani}}]\label{lem: BV lem}
Let $A=B(x_A,r_A)$ and $M=B(x_M,r_M)$ be arbitrary balls of $\R^k$ such that
\begin{equation*}
A\cap M\neq\varnothing\qquad\text{and}\qquad A\setminus(cM)\neq\varnothing
\end{equation*}
for some $c\geq3$. Then $r_M\leq r_A$ and $cM\subseteq5A$.
\end{lemma}
Consider any (sup norm) ball $B\subset\R^d$. If $B\cap K=\emptyset$ then trivially
\begin{equation*}
\mu(B)=0\leq|B|^{\Sigma-\varepsilon},
\end{equation*}
so we can suppose $B\cap K\neq\emptyset$. Let $n$ be the unique first level such that
\begin{equation*}
\#\{R\in\Rcal_{n-1}:R\cap B\neq\emptyset\}=1
\quad\text{and}\quad
\#\{R\in\Rcal_n:R\cap B\neq\emptyset\}\geq2.
\end{equation*}
If such level does not exist, then letting $R^{(n)}\in\Rcal_n$ denote the unique rectangle in the $n$th level that $B$ intersects, we have via Lemma~\ref{lem:rectangle-mass} that
\begin{equation*}
\mu(B)\leq\mu(R^{(n)})\leq C'|R^{(n)}|^{\Sigma-\varepsilon}\longrightarrow0
\qquad\text{as }n\to\infty,
\end{equation*}

so we are done. Therefore we may assume such level $n$ exists. Corresponding to level $n$ let $i=i(n)$ be the active coordinate and, for every $1\leq j\leq N$, let $B_j$ denote the projection of $B$ into the $j$th coordinate block. Consider the unique rectangle $R^{(n-1)}$ intersecting $B$ in $\Rcal_{n-1}$ and write $R^{(n-1)}=\prod_{j=1}^{N}I_j^{(n-1)}$.
If $|B|\geq|R^{(n-1)}|$ then, by Lemma~\ref{lem:rectangle-mass},
\begin{equation*}
\mu(B)\leq\mu(R^{(n-1)})\leq C'|R^{(n-1)}|^{\Sigma-\varepsilon}
\leq C'|B|^{\Sigma-\varepsilon},
\end{equation*}
and we are done. Henceforth suppose
\begin{equation}\label{eq: less than big rect}
|B|<|R^{(n-1)}|.
\end{equation}
Let
\begin{equation*}
F_i(B):=\{J_{i,a}\subset I_i^{(n-1)}:
\exists\,b=(b_j)_{j\neq i}\text{ such that }B\cap R_{a,b}\neq\emptyset\},
\end{equation*}
and, for each $J_{i,a}\in F_i(B)$ and $j\neq i$, let
\begin{equation*}
F_{j,a}(B):=\{J_{j,a,b_j}\in\Pcal_{j,a}:B_j\cap J_{j,a,b_j}\neq\emptyset\}.
\end{equation*}
These are the passive cubes occurring among the children with active index $a$ that intersect $B$.
\begin{itemize}
\item\textit{If $\#F_i(B)=1$,} say $J_{i,a}\in F_i(B)$, then by definition of level $n$
there is at least one passive coordinate $j\neq i$ with $\#F_{j,a}(B)\geq2$. That is, there is one `$a$' and several indices of `$b$'. In such $j$th coordinate let $A=B_j$ and $M=J_{j,a,b_j}\in F_{j,a}(B)$. Since the passive cubes are $3r$-pairwise disjoint, a point of $B_j$ in another cube lies outside $3M$. Lemma~\ref{lem: BV lem} is therefore applicable with $c=3$, and gives $|B_j|\geq|J_{j,a,b_j}|\asymp q_a^{-1-\tau_i}$.
By \eqref{eq:active-size} and \eqref{eq:passive-size}, all child cubes with this active index have side length at most $C|B|$, including those in passive coordinates where only one cube meets $B_j$.
Every such passive cube is consequently contained in a fixed dilation $CB_j$. Since the cubes in each passive family are disjoint, for every $j\neq i$ we have
\begin{equation}\label{eq: passive general ball cover}
\sum_{J_{j,a,b_j}\in F_{j,a}(B)}\Leb_{d_j}(J_{j,a,b_j})
\leq C|B|^{d_j}.
\end{equation}
So, using \eqref{eq:recursive-mass}, \eqref{eq:denom-bounds} and \eqref{eq:threshold},
\begin{align*}
\mu(B)
&\leq\mu(R^{(n-1)})
\frac{\Leb_{d_i}(B_i^+(\bp_a/q_a))}{\sum_a\Leb_{d_i}(B_i^+(\bp_a/q_a))}
\prod_{j\neq i}
\frac{\sum_{J_{j,a,b_j}\in F_{j,a}(B)}\Leb_{d_j}(J_{j,a,b_j})}
{\sum_{b_j}\Leb_{d_j}(J_{j,a,b_j})}\\
&\leq \kappa_{i} C^{N-1}\frac{\mu(R^{(n-1)})}{\Leb_d(R^{(n-1)})}
\left(q_a^{-1-\tau_i}\right)^{s_i}|B|^{\sum_{j\neq i}d_j}\\
&\leq \left(q_a^{-1-\tau_i}\right)^{s_i-\varepsilon}
|B|^{d-d_i}\\
&\leq C^{s_{i}-\varepsilon}|B|^{s_i+d-d_i-\varepsilon}
\\ &\leq C^{s_{i}-\varepsilon}|B|^{\Sigma-\varepsilon},
\end{align*}
where we have used $q_a^{-1-\tau_i}\leq C|B|$, $s_i-\varepsilon>0$ and $|B|<1$.

\item \textit{If $\#F_i(B)\geq2$,} then let $A=B_i$ and $M=B_i^-(\bp_a/q_a)\supset J_{i,a}\in F_i(B)$. We have chosen $G(R,n)$ large enough that
\[
c=\frac{\rho_i(q_a)}{r_i(q_a)}\geq3,
\qquad cM=B_i^+(\bp_a/q_a).
\]
By Lemma~\ref{lem:filling}, the balls $B_i^+(\bp_a/q_a)$ are disjoint over $a$. Since $B_i$ meets an active small cube with a different index, $B_i\setminus(cM)\neq\emptyset$.
Lemma~\ref{lem: BV lem}, with this value of $c$, gives
\[
B_i^+(\bp_a/q_a)\subseteq5B_i,
\qquad q_a^{-1-\tau_i}\leq\tfrac12|B|.
\]
Consequently
\begin{equation*}
\sum_{J_{i,a}\in F_i(B)}\Leb_{d_i}(B_i^+(\bp_a/q_a))
\leq\Leb_{d_i}(5B_i).
\end{equation*}
Again we have \eqref{eq: passive general ball cover} for every $a$ and every passive coordinate, even when $\#F_{j,a}(B)=1$.
Therefore, using \eqref{eq:recursive-mass} we have
\begin{align*}
\mu(B)
&\leq\frac{\mu(R^{(n-1)})}{\sum_a\Leb_{d_i}(B_i^+(\bp_a/q_a))}
\sum_{J_{i,a}\in F_i(B)}\Biggl[\Leb_{d_i}(B_i^+(\bp_a/q_a))\times\prod_{j\neq i}
\frac{\sum_{J_{j,a,b_j}\in F_{j,a}(B)}\Leb_{d_j}(J_{j,a,b_j})}
{\sum_{b_j}\Leb_{d_j}(J_{j,a,b_j})}\Biggr]\\
&\leq C\frac{\mu(R^{(n-1)})}{\Leb_d(R^{(n-1)})}|B|^d\\
&\leq C|R^{(n-1)}|^{\Sigma-\varepsilon-d}|B|^d\\
&=C|B|^{\Sigma-\varepsilon}
\left(\frac{|B|}{|R^{(n-1)}|}\right)^{d-\Sigma+\varepsilon}\\
&\leq C|B|^{\Sigma-\varepsilon}.
\end{align*}
Here we used Lemma~\ref{lem:rectangle-mass} and \eqref{eq:shape}, which gives $\Leb_d(R^{(n-1)})\asymp|R^{(n-1)}|^d$.
The last inequality follows from \eqref{eq: less than big rect} and $\Sigma<d$.
\end{itemize}
This completes all cases. Taking $C'>0$ to be the maximum constant appearing in each of the above cases, we have that for any sufficiently small ball $B\subset\R^d$ 
\begin{equation*}
\mu(B)\leq C'|B|^{\Sigma-\varepsilon},
\end{equation*}
and so, by the mass distribution principle,
\[
\dimH\prod_{i=1}^{N}\Bad_{d_i}(\tau_i)\geq\dimH K\geq\Sigma-\varepsilon.
\]
Since $\varepsilon>0$ was arbitrary the proof is complete.

\end{document}